\documentclass[a4paper,10pt]{article}
\usepackage{amsfonts}
\usepackage{amssymb,amsthm}
\usepackage{amsmath}

\newcommand{\R}{\mathbb{R}}
\newcommand{\N}{\mathbb{N}}

\newcommand{\Z}{\mathbb{Z}}

\newtheorem{theorem}{Theorem}
\newtheorem{lemma}[theorem]{Lemma}

\newcommand{\be}{\begin{equation}}
\newcommand{\ee}{\end{equation}}
\newcommand{\bee}{\begin{equation*}}
\newcommand{\eee}{\end{equation*}}
\newcommand{\bp}{\begin{proof}}
\newcommand{\ep}{\end{proof}}

\def\tbigcup{\bigcup}
\begin{document}
\title{Hardy inequality and heat semigroup estimates for Riemannian manifolds with singular
data}
\author{{M. van den Berg, P. Gilkey\thanks{Partially supported by Project MTM2009-07756 (Spain)},
 A. Grigor'yan\thanks{ Partially supported by SFB701 (Germany)}, K. Kirsten
\thanks{Supported by National Science Foundation Grant PHY-0757791}} \\
\\
School of Mathematics, University of Bristol\\
University Walk, Bristol BS8 1TW, UK\\
\texttt{M.vandenBerg@bris.ac.uk}\\ \\
Mathematics Department, University of Oregon\\
Eugene, OR 97403, USA\\
\texttt{gilkey@uoregon.edu}\\ \\
Fakult\"at f\"ur Mathematik, Universit\"at Bielefeld\\
Postfach 100131, D-33501 Bielefeld, Germany\\
\texttt{grigor@math.uni-bielefeld.de}\\ \\
Department of Mathematics, Baylor University\\
Waco, Texas, TX 76798, USA\\
\texttt{Klaus\_Kirsten@baylor.edu}}
\date{2 June 2011}\maketitle
%\vskip 3truecm \indent
\begin{abstract}\noindent Upper bounds are obtained for the
heat content of an open set $D$ in a geodesically complete
Riemannian manifold $M$ with Dirichlet boundary condition on
$\partial D$, and non-negative initial condition. We show that
these upper bounds are close to being sharp if (i) the
Dirichlet-Laplace-Beltrami operator acting in $L^2(D)$ satisfies a
strong Hardy inequality with weight $\delta^2$, (ii) the initial
temperature distribution, and the specific heat of $D$ are given
by $\delta^{-\alpha}$ and $\delta^{-\beta}$ respectively, where
$\delta$ is the distance to $\partial D$, and $1<\alpha<~2,
1<\beta <~2$.
\end{abstract}
\vskip 1truecm \noindent \ \ \ \ \ \ \ \  { Mathematics Subject
Classification (2000)}: 58J32; 58J35; 35K20.
\begin{center} \textbf{Keywords}: Hardy inequality, heat content, singular data.\\
\end{center} \mbox{}\newpage
\section{Introduction}
\label{sec1}Let $D$ be a smooth, connected, $m$- dimensional
Riemannian manifold and let $\Delta $ be the Laplace-Beltrami
operator on $D$. It is well known (see \cite{EBD3}, \cite{GB})
that the heat equation
\begin{equation}
\Delta u=\frac{\partial u}{\partial t},\quad x\in D,\quad t>0,
\label{e1}
\end{equation}%
has a unique minimal positive fundamental solution $p(x,y;t)$
where $ x\in D$, $y \in D$, $t>0$. This solution, the Dirichlet
heat kernel for $D$, is symmetric in $x,y$, strictly positive,
jointly smooth in $x,y\in D$ and $t>0$, and it satisfies the
semigroup property
\begin{equation}
p(x,y;s+t)=\int_{D}\ p(x,z;s)p(z,y;t)dz , \label{e4}
\end{equation}%
for all $x,y\in D$ and $t,s>0$, where $dz $ is the Riemannian
measure on $D$. Equation (\ref{e1}) with the initial condition
\begin{equation}
u(x;0^+)=\psi (x),\quad x\in D,  \label{e2}
\end{equation}%
has a solution
\begin{equation}
u_{\psi }(x;t)=\int_{D}\ p(x,y;t)\psi(y)dy  , \label{e5}
\end{equation}
for any function $\psi $ on $D$ from a variety of function spaces
like $C_{b}\left( D\right)$ or  $ L^{p}\left( D\right) $, $1\leq
p< \infty $. Note that $u_{\psi}\in C_{b}\left( D\right) $ if
$\psi \in C_{b}\left( D\right) $ or that $u_{\psi}\in L^{p}\left(
D\right)$ if $\psi \in L^{p}\left( D\right)$. Initial condition
(\ref{e2}) is understood in the sense that $u_{\psi}\left(
\cdot;t\right) \rightarrow \psi \left( \cdot\right) $ as
$t\rightarrow 0^+$, where the convergence is appropriate for the
function space of initial conditions. For example, if $\psi \in
C_{b}\left( D\right) $ then the convergence is locally uniform, or
if $\psi \in L^{p}\left( D\right) $, $1\leq p< \infty $ then the
convergence is in the norm of $L^{p}\left( D\right) $. In general,
(\ref{e5}) is not the unique solution of (\ref{e1})-( \ref{e2}).
However, it has the following distinguished property: if $\psi
\geq 0$ then $u_{\psi }$ is the minimal non-negative solution of
that problem (and if $\psi $ is signed then $u_{\psi }=u_{\psi
_{+}}-u_{\psi _{-}} $). If $D$ is an open subset of another
Riemannian manifold $M$ and if the boundary $\partial D$ of $D$ in
$M$ is smooth then the minimality property of $u_{\psi }$ implies
that, for any $t>0$,
\begin{equation}
\lim_{x\rightarrow \partial D}u_{\psi }(x;t)=0.\   \label{e3}
\end{equation}%
If $\partial D$ is non-smooth then (\ref{e3}) can still be
understood in a weak sense. Expression (\ref{e5}) makes sense for
any non-negative measurable function $\psi$ on $D$, provided the
value $+\infty $ is allowed for $u_{\psi }$. It is known that if
$u_{\psi }\in L_{loc}^{1}\left( D\times \mathbb{R}_{+}\right) $
then $u_{\psi }$ is a smooth function in $ D\times \mathbb{R}_{+}$
and it solves (\ref{e1}) (see p. 201 in \cite{GB}). For any two
non-negative measurable functions $\psi _{1},\psi _{2}$ on $D$, we
define for $t>0$
\begin{equation}
Q_{\psi _{1},\psi _{2}}(t)=\iint_{D\times D}
p(x,y;t)\psi_1(x)\psi_2(y)dxdy. \label{e6}
\end{equation}%
Using the properties of the Dirichlet heat kernel we have for
$0<s<t$
\begin{equation} \label{12}
Q_{\psi _{1},\psi _{2}}(t)=\int_{D}\ u_{\psi _{1}}(x;s) u_{\psi
_{2}}(x;t-s)dx.
\end{equation}
Assuming that $D$ is an open subset of a complete Riemannian
manifold $M$, $Q_{\psi _{1},\psi _{2}}(t)$ has the following
physical interpretation: it is the amount of heat in $D$ at time
$t$ if $D$ has initial temperature distribution $\psi _{1}$, and a
specific heat $\psi _{2}$, while the $\partial D$ is kept at fixed
temperature $0$.

This function has been subject of a thorough investigation. Its
asymptotic behavior for small $t$ is well understood if $D$ has
compact closure with $C^{\infty }$ boundary, and both $\psi _{1}$
and $\psi _{2}$ are $C^{\infty }$ on the closure $\overline{D}$ of
$D$. In that case $ Q_{\psi _{1},\psi _{2}}(t)$ has an asymptotic
series in $t^{1/2}$, and its coefficients are computable in terms
of local geometric invariants \cite {vdB1,PG}. No such series are
known if $D$ is unbounded, or if either the initial data or
$\partial D$ are non-smooth.

In this paper we will obtain upper bounds for the heat content
$Q_{\psi _{1},\psi _{2}}(t)$ under quite general assumptions on
$D$ and on $\psi _{1}$ and $\psi _{2}$.

We are particularly interested in the situation where $D$ is a
open subset of another manifold $M$, and  where $\psi _{1}\left(
x\right) $ and $\psi _{2}\left( x\right) $ blow up as
$x\rightarrow \partial D$. In
order to guarantee finite heat content for $t>0$, sufficient cooling at $%
\partial D$ needs to take place. This will be guaranteed by a
condition on $D$, that is formulated in terms of a Hardy
inequality. Note that in this setting $Q_{\psi _{1},\psi
_{2}}\left( t\right) $ may be unbounded as $t\rightarrow 0^+$, and
one of the interesting points of this study is to obtain the rate
of convergence of $Q_{\psi _{1},\psi _{2}}\left( t\right) $ to
$+\infty $ as $t\rightarrow 0^+.$

Given a positive measurable function $h$ on a manifold $D$, we say
that the Dirichlet Laplacian acting in $L^2(D)$ satisfies a strong
Hardy inequality with weight $h$ if, for all $w\in C_{c}^{\infty
}(D)$,
\begin{equation}
\int_{D}|\nabla w|^{2}\geq \int_{D}\ \frac{w^{2}}{h}\ .
\label{e7}
\end{equation}%
Here, and in what follows, we put $\int_Df=\int_Df(x)dx$ if this
does not cause confusion. We also put $|D|=\int_D1$, and $\Vert
f\Vert_p=\left(\int_D|f|^p\right)^{1/p}$. A typical example of a
Hardy inequality
is when $D$ is an open subset of another manifold $M$, and%
\begin{equation}
h(x)=c^{2}\delta (x)^{2},  \label{e8}
\end{equation}%
where $c\geq 2$ is a constant, $\delta$ is the distance to the
boundary,
\begin{equation*}
\delta (x)=\min \{d(x,y):y\in \partial D\},  \label{e9}
\end{equation*}
and $d(x,y)$ is the geodesic distance from $x$ to $y$ on $M$. Both
the validity and applications of Hardy inequalities with weight
(\ref{e8}) have been investigated extensively \cite{A}, \cite{B},
\cite{EBD1}, \cite{EBD2}, \cite{EBD3}, \cite{vdB3}. For example,
inequality \eqref{e7} holds with weight (\ref{e8}) with $c=4$ if
$D$ is simply connected with
non-empty boundary in $\mathbb{R}^{2}$, with $c=2$ if $D$ is convex in $%
\mathbb{R}^{m}$, and for some $c\geq 2$ if $D$ is bounded with
smooth boundary in $\mathbb{R}^{m}$.

In \cite{vdB2} it was shown that if $D$ has finite measure and
satisfies the Hardy inequality with weight $h$, and if $\psi $ is
a non-negative measurable function on $D$, such that, for some
$q>1$,
\begin{equation}
\lVert \psi h^{1/q}\rVert _{q/(q-1)}<\infty ,  \label{e10}
\end{equation}%
then, for all $t>0,$
\begin{equation}\label{e11}
Q_{\psi ,1}(t)\leq \left( \frac{q^{2}}{4(q-1)}\right) ^{1/q}\lVert
\psi h^{1/q}\rVert _{q/(q-1)}\lVert 1-u_{1}\left( \cdot ;t\right)
\rVert _{1}^{1/q}t^{-1/q},
\end{equation}
where $Q_{1,1}$ is defined by (\ref{e6}) for $\psi _{1}=\psi
_{2}=1$, that is,
\begin{equation*}\label{e12}
Q_{1,1}(t)=\int_{D}u_{1}\left( x;t\right) dx=\iint_{D\times
D}p(x,y;t)dxdy.
\end{equation*}

A similar estimate holds for arbitrary open sets $D\subset
\mathbb{R}^{m},$ satisfying the Hardy inequality with weight $h$.
If $\psi $ is a non-negative measurable
function on $D$ such that, for some $q>1$,%
\begin{equation}
\lVert \max\{\psi, 1\}h^{1/q}\rVert _{q/\left( q-1\right) }<\infty
, \label{e13}
\end{equation}
then, for all $t>0$,
\begin{equation}
Q_{\psi ,1}(t)\leq a(q)\lVert \psi h^{1/q}\rVert _{q/(q-1)}\lVert
h^{1/(q(q-1))}\rVert _qt^{-1/(q-1)}, \label{e14}
\end{equation}
where
\begin{equation}
a(q)=4^{-1/q}\left( \frac{q}{q-1}\right) ^{(2q-1)/(q(q-1))}.
\label{e15}
\end{equation}

Below we give a sufficient condition for the finiteness of
$Q_{\psi _{1},\psi _{2}}(t)$ for all $t>0$, and reduce the problem
of finding upper bounds for $Q_{\psi _{1},\psi _{2}}(t)$ to the
case $\psi _{1}=\psi _{2}$.
\begin{theorem}\label{The0}
Let $\psi_1$ and $\psi_2$ be non-negative and Borel measurable on
a manifold $D$.
\begin{enumerate} \item[(i)] If $Q_{\psi_i,\psi_i}(t)<\infty,
i=1,2,$ for all $t>0$, then $Q_{\psi_1,\psi_2}(t)<\infty$ for all
$t>0$, and
\begin{equation}\label{e17}
Q_{\psi_1,\psi_2}(t)\le \left(
Q_{\psi_1,\psi_1}(t)Q_{\psi_2,\psi_2}(t)\right)^{1/2}, \ t>0 .
\end{equation}
\item[(ii)]If $Q_{\psi_i,1}(t)<\infty, i=1,2,$ for all $t>0$,
and if \begin{equation}\label{e17a}c_t:= \sup_{x \in D}p(x,x;t)<
\infty, \ \ t>0,\end{equation}
 then
\begin{equation*}\label{e16}
Q_{\psi_1,\psi_2}(t)\le c_{t/3}Q_{\psi_1,1}(t/3)Q_{\psi_2,1}(t/3)<
\infty, \ t>0.
\end{equation*}
\end{enumerate}
\end{theorem}

Our main results are the following three theorems, in which we
assume that $D$ is a Riemannian manifold that satisfies the Hardy
inequality with some weight $h$, and $\psi $ is a non-negative
measurable function on $D$. In particular we do not assume any
smoothness conditions on $\partial D$, nor do we assume that $D$
has finite measure or that $D$ is bounded.
\begin{theorem}
\label{The1}If $|D| <\infty $, and if there exists $%
1<q\leq 2$ such that
\begin{equation}
\lVert \psi h^{1/q}\rVert _{q/(q-1)}<\infty,  \label{e18}
\end{equation}
then, for all $t>0$,
\begin{equation}
Q_{\psi ,\psi }(t)\leq \frac{q^{(4-q)/q}}{(2(q-1))^{2/q}}\lVert
\psi h^{1/q}\rVert _{q/(q-1)}^{2}\lVert 1-u_{1}\left( \cdot
;t\right) \rVert _{1}^{(2-q)/q}t^{-2/q}.  \label{e19}
\end{equation}
\end{theorem}

\begin{theorem}
\label{The2}If $1<q\leq 2$ is such that \emph{(\ref{e18})} holds
and that
\begin{equation*}
\lVert h^{1/q}\rVert _{q/(q-1)}<\infty ,
\end{equation*}
then
\begin{equation}
Q_{\psi ,\psi }(t)\leq b(q)\ \lVert \psi h^{1/q}\rVert
_{q/(q-1)}^{2}\ \lVert h^{1/q}\rVert
_{q/(q-1)}^{(2-q)/(q-1)}t^{-1/(q-1)},\ t>0,  \label{e20}
\end{equation}%
where
\begin{equation}
b(q)=2^{(4-3q)/(q(q-1))}\left( \frac{q}{q-1}\right)
^{(q^{2}-4q+2)/(q(1-q))}. \label{e21}
\end{equation}
\end{theorem}

\begin{theorem}
\label{The3}If $0\leq r\leq 2$, and $1<q\leq 2$ are such that
\begin{equation*}
||\psi ^{r}||_{q}<\infty,
\end{equation*}
and
\begin{equation*}
\lVert \psi ^{2-r}h^{1/q}\rVert _{(q-1)/q}<\infty ,  \label{e22}
\end{equation*}
then
\begin{equation}
Q_{\psi ,\psi }(t)\leq \left( \frac{q}{4(q-1)}\right) ^{1/q}\lVert
\psi ^{r}\rVert _{q}\lVert \psi ^{2-r}h^{1/q}\rVert
_{q/(q-1)}t^{-1/q},\ t>0. \label{e23}
\end{equation}
\end{theorem}

In Theorem \ref{The5} in Section \ref{sec3} we use the bounds of
Theorems \ref{The1} and \ref{The3} together with (\ref{e17}) to
obtain an upper bound for the heat content of $D$, when $D$
satisfies a Hardy inequality with weight (\ref{e8}), and $\psi
_{1}(x)=\delta(x) ^{-\alpha }$ and $\psi _{2}(x)=\delta(x)
^{-\beta }$, where $\alpha ,\beta \in \left( 1,2\right) $. Even
though the bounds in e.g. \ref{The1} and \ref{The3} look very
different, both of them are needed to cover the maximal range of
$\alpha$ and $\beta$ in Theorem \ref{The5}.

Theorem \ref{The1} has a curious consequence. We claim that if a
manifold $D$ has finite measure $|D|$, and is stochastically
complete then no Hardy inequality holds on $D$ (which confirms the
philosophy that the Hardy inequality corresponds to cooling that
comes from the boundary). Indeed, stochastic completeness means
that $u_{1}\equiv 1$. In this
case, $\left\Vert 1-u_1(\cdot;t)\right\Vert _{1}=0$ so that we obtain from (\ref%
{e19}) that $Q_{\psi ,\psi }\left( t\right) =0$ whenever function
$\psi $ satisfies the condition (\ref{e18}) for some $q\in \left(
1,2\right) $. However, if $h$ is finite then it is easy to
construct a non-trivial function $\psi $ that satisfies
(\ref{e18}): choose any
measurable set $S$ with finite positive measure such that $h$ is bounded on $%
S$, and let $\psi =1_{S}$. Then (\ref{e18}) holds with any $q>1$ while $%
Q_{\psi ,\psi }\left( t\right) >0$ so that we obtain
contradiction. Of course, without the finiteness of $|D|$, the
Hardy inequality may hold on stochastically complete manifolds like $\mathbb{R}%
^{m}\setminus \{0\}$.

This paper is organized as follows. In Section \ref{sec2} we will
prove Theorems \ref{The0}, \ref{The1}, \ref{The2} and \ref{The3}.
In Section \ref{sec3} we will state and prove Theorem \ref{The5}.
Finally in Section \ref{sec4} we obtain very refined asymptotics
in the special case of the ball in $\R^3$ with $\psi
_{1}(x)=\delta(x) ^{-\alpha }, \alpha<2, \psi _{2}(x)=\delta(x)
^{-\beta },\beta<2$, and $\alpha + \beta
>3$ (Theorem \ref{The6}). This special case shows that the bound obtained in
Theorem \ref{The5} is close to being sharp. Moreover it suggests
formulae for the first few terms in the asymptotic series of a
compact Riemannian manifold $D$ with the singular data above.

\section{Proofs of Theorems \ref{The0}, \ref{The1}, \ref{The2} and \ref{The3}}

\label{sec2}

\begin{proof}[Proof of Theorem \protect\ref{The0}]
In both parts, it suffices to prove the claims for non-negative functions $%
\psi _{1},\psi _{2}$ from $L^{2}\left( D\right) $. Arbitrary
non-negative measurable functions $\psi _{1},\psi _{2}$ can be
approximated by monotone increasing sequences of non-negative
functions from $L^{2}\left( D\right) $, whence the both claims
follow by the monotone convergence theorem.

To prove part (i) we use symmetry and the semigroup property, and
obtain by \eqref{12} for $s=t/2$ that
\begin{align*}
Q_{\psi _{1},\psi _{2}}(t)&=\int_Du_{\psi_1}(x;t/2)u_{\psi_2}(x;t/2)dx\nonumber \\
&\le\left(\int_Du_{\psi_1}^2(x;t/2)dx\right)^{1/2}\left(\int_Du_{\psi_2}^2(x;t/2)dx\right)^{1/2}\nonumber
\\ &=\left(Q_{\psi _{1},\psi _{1}}(t)Q_{\psi _{2},\psi
_{2}}(t)\right)^{1/2}.
\end{align*}
It follows from (\ref{e4}) and \eqref{e17a} that%
\begin{equation}
p(x,y;t)\leq (p(x,x;t)p(y,y;t))^{1/2}\leq c_{t}.  \label{e23c}
\end{equation}%

To prove part (ii) we have by \eqref{e23c} that
\begin{align}
p(x,y;t)&=\iint_{D\times
D}p(x,z_1;t/3)p(z_1,z_2;t/3)p(z_2,y;t/3)dz_1dz_2\nonumber
\\ &\le c_{t/3}u_1(x;t/3)u_1(y;t/3).
\end{align}
This together with definition \eqref{e6} completes the proof.
\end{proof}

For the proofs of Theorems \ref{The1}, \ref{The2}, \ref{The3},
choose a sequence $\left\{ D_{n}\right\} $ that consists of
precompact
open subsets of $D$ with smooth boundaries such that $\overline{D}%
_{n}\subset D_{n+1}$ and $\tbigcup_{n}D_{n}=D.$ Obviously, Hardy
inequality (\ref{e7}) remains true in any $D_{n}$ with the same
weight $h$, because $C_{c}^{\infty }\left( D_{n}\right) \subset
C_{c}^{\infty }\left( D\right) .$ Moreover, we claim that
(\ref{e7}) holds for any function $w\in C\left(
\overline{D}_{n}\right)\cap C^{1}\left( D_{n}\right) $ that
satisfies the boundary condition $w|_{\partial D_{n}}=0.$ Indeed, if $%
\int_{D_{n}}\left\vert \nabla w\right\vert ^{2}=\infty $ then
(\ref{e7}) is trivially satisfied. If $\int_{D_{n}}\left\vert
\nabla w\right\vert ^{2}<\infty $ then $w$ belongs to the Sobolev
space $W^{1,2}\left(
D_{n}\right) $. Extend function $w$ to $D_{n+1}$ by setting $w=0$ in $%
D_{n+1}\setminus \overline{D}_{n}$. Due to the boundary condition $%
w|_{\partial D_{n}}=0$, we obtain that $w_{n}\in W^{1,2}\left(
D_{n+1}\right) $. Since $w$ is compactly supported in $D_{n+1}$,
it follows that $w\in W_{0}^{1,2}\left( D_{n+1}\right) $ where
$W_{0}^{1,2}\left( \Omega \right) $ is the closure $C_{c}^{\infty
}\left( \Omega \right) $ in $W^{1,2}\left( \Omega \right) $. Since
the Hardy inequality (\ref{e7}) holds for functions from
$C_{c}^{\infty}(D_{n+1})$, passing to the limit in
$W^{1,2}(D_{n+1})$ and using Fatou's lemma, we obtain that $w$
also satisfies (\ref{e7}).

Assume for a moment that the statements of the theorems have been
proved in each domain $D_{n}$. Then one can take the limit in
(\ref{e19}), (\ref {e20}), (\ref{e23}) as $n\rightarrow \infty $,
and obtain the statements for $ D$. Indeed, the left hand side of
these inequalities is $Q^{D_n}_{\psi ,\psi }\left(
t\right)=\iint_{D_n\times D_n} p_{D_n}(x,y;t)\psi(x)\psi(y)dxdy$,
where $p_{D_n}$ is the Dirichlet heat kernel for $D_n$. This
converges to $Q_{\psi ,\psi }^{D}\left( t\right) $ as
$n\rightarrow \infty $. The right hand sides of (\ref{e19}),
(\ref{e20}), (\ref{e23}) contain various $L^{p}\left(
D_{n}\right) $-norms that can be estimated from above by the $%
L^{p}\left( D\right) $-norms. The only exception is the term
$\lVert 1-\int_{D_n}p_{D_n}(\cdot,y;t)dy\rVert _{1}$ in
(\ref{e19}) that is decreasing as $n\rightarrow \infty $. If $|D|
<\infty $ then $1\in L^{1}\left( D\right) $ so that the passage to
the limit is justified by the dominated convergence theorem.

Hence, it suffices to prove each of the statements for $D_{n}$
instead of $D$. Renaming $D_{n}$ back to $D$, we assume in all
three proofs that $D$ is a precompact open domain with smooth
boundary in $M$.

Another observation is that all inequalities (\ref{e19}), (\ref{e20}), (\ref%
{e23}) survive the increasing monotone limits in $\psi $. So it
suffices to prove them when $\psi $ is bounded and has a compact
support in $D$, which
will be assumed below. Furthermore, since all the statements of Theorems \ref%
{The1}, \ref{The2}, \ref{The3} are homogeneous with respect to
$\psi $, we can assume that $0\leq \psi \leq 1$. If $\psi \equiv
0$ then there is nothing to prove; hence, we assume that $\psi $
is non-trivial. Then $u_{\psi }\left( x;t\right) $ is smooth and bounded in $\overline{D}%
\times \left( 0,+\infty \right) $ and positive in $D\times \left(
0,+\infty \right) $.

\begin{proof}[Proof of Theorem \protect\ref{The1}]
Let $\nu $ be the outwards normal vector field on $\partial D$.
Using the Green's formula, we obtain
\begin{eqnarray}\label{L4}
-\frac{d}{dt}\int_{D}u_{\psi}^{q}
&=&-q\int_{D}u_{\psi}^{q-1}\frac{\partial u_{\psi}}{\partial t}
\notag \\
&=&-q\int_{D}u_{\psi}^{q-1}\Delta u_{\psi}  \notag \\
&=&-q\int_{\partial D}u_{\psi}^{q-1}\frac{\partial u_{\psi}}{\partial \nu }%
+q\int_{D}\left( \nabla u_{\psi}^{q-1},\nabla u_{\psi}\right)  \notag \\
&=&q\left( q-1\right) \int_{D}u_{\psi}^{q-2}\left\vert \nabla
u_{\psi}\right\vert ^{2}, %
\end{eqnarray}
where we have used that $q>1$ and, hence $u_{\psi}^{q-1}=0$ on
$\partial D$. Observing that $u_{\psi}^{q/2}\in C\left(
\overline{D}\right) \cap C^{1}\left( D\right) $,
\begin{equation*}
\left\vert \nabla u_{\psi}^{q/2}\right\vert
^{2}=\frac{q^{2}}{4}u_{\psi}^{q-2}\left\vert \nabla
u_{\psi}\right\vert ^{2},
\end{equation*}
and applying the Hardy inequality (\ref{e7}) to $u^{q/2}$, we
obtain that
\begin{equation}
-\frac{d}{dt}\int_{D}u_{\psi}^{q}=\frac{4(q-1)}{q}\int_{D}|\nabla
(u_{\psi}^{q/2})|^{2}\geq
\frac{4(q-1)}{q}\int_{D}\frac{u_{\psi}^{q}}{h}. \label{e24}
\end{equation}%
By H\"{o}lder's inequality we have that
\begin{eqnarray}
Q_{\psi ,\psi }(t) &=&\int_{D}u_{\psi}\psi  \notag \\
&\leq &\left( \int_{D}\left( \frac{u_{\psi}}{h^{1/q}}\right)
^{q}\right)
^{1/q}\left( \int \left( \psi h^{1/q}\right) ^{\frac{q}{q-1}}\right) ^{\frac{%
q-1}{q}}  \notag \\
&=&\left( \int_{D}\frac{u_{\psi}^{q}}{h}\right) ^{1/q}\lVert \psi
h^{1/q}\rVert _{q/(q-1)}.  \label{e25}
\end{eqnarray}
By (\ref{e24}) and (\ref{e25}) we conclude that
\begin{equation}
-\frac{d}{dt}\int_{D}u_{\psi}^{q}\geq \frac{4(q-1)}{q}\lVert \psi
h^{1/q}\rVert _{q/(q-1)}^{-q}\left( Q_{\psi ,\psi }(t)\right)
^{q}.  \label{e26}
\end{equation}
Note that the function $t\mapsto Q_{\psi ,\psi }(t)=\left\Vert
u_{\psi}\left( \cdot ;t/2\right) \right\Vert ^2_2$ is decreasing
in $t$, which, for example,
follows from (\ref{L4}) with $q=2$. Integrating differential inequality (%
\ref{e26}) with respect to $t$ over the interval $[t,2t]$ gives
that
\begin{equation}
\int_{D}u_{\psi}^{q}\geq \frac{4(q-1)}{q}\lVert \psi h^{1/q}\rVert
_{q/(q-1)}^{-q}\left( Q_{\psi ,\psi }(2t)\right) ^{q}t.
\label{e27}
\end{equation}

On the other hand, using $1<q\leq 2$ and the H\"{o}lder inequality, we obtain%
\begin{equation*}
\int_{D}u_{\psi}^{q}=\int_{D}u_{\psi}^{2-q}u_{\psi}^{2q-2}\leq
\left( \int_{D}u_{\psi}\right) ^{2-q}\left(
\int_{D}u_{\psi}^{2}\right) ^{q-1}
\end{equation*}
that is,%
\begin{equation}
\int_{D}u_{\psi}^{q}\leq \left( Q_{\psi ,1}\left( t\right) \right)
^{2-q}\left( Q_{\psi ,\psi }(2t)\right) ^{q-1}.  \label{e28}
\end{equation}%
Combining (\ref{e27}) and (\ref{e28}) yields%
\begin{equation}
Q_{\psi ,\psi }(2t)\leq \frac{q}{4(q-1)}||\psi
h^{1/q}||_{q/(q-1)}^{q}\left( Q_{\psi ,1}\left( t\right) \right)
^{2-q}t^{-1}.  \label{e29}
\end{equation}
Estimating $Q_{\psi ,1}$ by \eqref{e11}, we obtain%
\begin{equation*}
Q_{\psi ,\psi }(2t)\leq \frac{q^{(4-q)/q}}{(4(q-1))^{2/q}}\lVert
\psi h^{1/q}\rVert _{q/(q-1)}^{2}\lVert 1-u_{1}(\cdot;t)\rVert
_{1}^{(2-q)/q}t^{-2/q}, \label{e30}
\end{equation*}
which completes the proof.
\end{proof}

\begin{proof}[Proof of Theorem \protect\ref{The2}]
Since $\psi \leq 1$ we have that \eqref{e13} is satisfied. We
obtain by \eqref{e14} and \eqref{e29} that
\begin{equation*}
Q_{\psi ,\psi }(2t)\leq \frac{q}{4(q-1)}a(q)^{2-q}\lVert \psi
h^{1/q}\rVert _{q/(q-1)}^{2}\ \lVert h^{1/q}\rVert
_{q/(q-1)}^{(2-q)/(q-1)}t^{-1/(q-1)}. \label{e33}
\end{equation*}
This completes the proof of Theorem \ref{The2} since, by
(\ref{e15}) and (\ref{e21}),
\begin{equation*}
2^{1/(q-1)}\frac{q}{4(q-1)}a(q)^{2-q}=b(q).  \label{e34}
\end{equation*}
\end{proof}

\begin{proof}[Proof of Theorem \protect\ref{The3}]
By the arithmetic-geometric mean inequality, we have%
\begin{equation*}
\psi (x)\psi (y)\leq \frac{1}{2}\left( \psi (x)^{r}\psi
(y)^{2-r}+\psi (x)^{2-r}\psi (y)^{r}\right).  \label{e35}
\end{equation*}
By non-negativity and symmetry of the Dirichlet heat kernel
\begin{equation}
Q_{\psi ,\psi }\left( t\right) \leq \int_{D}u_{\psi ^{r}}\psi
^{2-r}.  \label{e36}
\end{equation}%
Next, H\"{o}lder's inequality yields%
\begin{equation}
\int_{D}u_{\psi ^{r}}\psi ^{2-r}\leq \left( \int_{D}u_{\psi ^{r}}^{q}\frac{1%
}{h}\right) ^{1/q}\lVert \psi ^{2-r}h^{1/q}\rVert _{q/(q-1)}.
\label{e37}
\end{equation}%
By (\ref{e24}) we have
\begin{equation}
-\frac{d}{dt}\int_{D}u_{\psi ^{r}}^{q}\geq
\frac{4(q-1)}{q}\int_{D}u_{\psi ^{r}}^{q}\frac{1}{h}.  \label{e38}
\end{equation}%
Combining (\ref{e36}), (\ref{e37}), (\ref{e38}) we obtain that
\begin{equation*}
\left( Q_{\psi ,\psi }\left( t\right) \right) ^{q}\leq -\frac{q}{4(q-1)}%
\frac{d}{dt}\left( \int_{D}u_{\psi ^{r}}^{q}\right) \lVert \psi
^{2-r}h^{1/q}\rVert _{q/(q-1)}^{q}.  \label{e39}
\end{equation*}
Since the function $t\mapsto Q_{\psi ,\psi }\left( t\right) $ is
decreasing in $t$, we obtain by integrating the differential
inequality (\ref{e38}) with respect to $t$ over the interval
$[0,t]$ that
\begin{equation*}
t\left( Q_{\psi ,\psi }\left( t\right) \right) ^{q}\leq \frac{q}{4(q-1)}%
\left( \int_{D}\psi ^{rq}\right) \lVert \psi ^{2-r}h^{1/q}\rVert
_{q/(q-1)}^{q},  \label{e40}
\end{equation*}
and (\ref{e23}) follows.
\end{proof}
\section{Singular initial temperature and singular specific heat \label{sec3}}
Below we make some further hypothesis on the geometry of $D$, and
obtain an upper bound for the heat content for a wide class of
geometries using Theorems \ref{The1} and \ref{The3}, and
\eqref{e17}, if the initial temperature distribution and specific
heat are given by $\delta^{-\alpha}, 1<\alpha<2$, and
$\delta^{-\beta}, 1< \beta<2$ respectively.
\begin{theorem}\label{The5}Let $D$ be an open set in a smooth complete
$m$-dimensional Riemannian manifold $M$, and suppose that
\begin{enumerate}
\item[i.] The Ricci curvature on $M$ is non-negative.
\item[ii.]For $x\in D$, \begin{equation*}\label{e41}
\psi_{\alpha}(x)=\delta(x)^{-\alpha}.
\end{equation*}
\item[iii.]$D$ has finite inradius i.e.
$\rho_D=\sup\{\delta(x): x\in D\}< \infty$.
\item[iv.]There
exist constants $\kappa_D<\infty, d\in [m-1,m)$ such that
\begin{equation}\label{e42}
\int_{\{ x \in D:\delta(x)<\rho\}}1 \le \kappa_D \rho^{m-d},\
0<\rho \le \rho_D.
\end{equation}
\item[v.]The strong Hardy inequality \eqref{e7} holds with \eqref{e8} for some
$c\ge 2$.
\end{enumerate}
If $1<\alpha<2, 1<\beta<2$, and if $\epsilon>0$ then
\begin{equation}\label{e43}
Q_{\psi_{\alpha},\psi_{\beta}}(t)=O(t^{-\epsilon+(m-d-\alpha-\beta)/2}),\
t\rightarrow 0 .
\end{equation}
\end{theorem}
\begin{proof}
Note that (iii) and (iv) in Theorem \ref{The5} imply that $|D|\le
\kappa_D\rho_D^{m-d}<\infty.$ By \eqref{e17} it suffices to prove
\eqref{e43} in the special case $\alpha=\beta$ with $1<\alpha<2$.
In order to estimate $\lVert 1-u_1(\cdot;t)\rVert_1$ in Theorem
\ref{The1} we rely on the following lower bound for $u_1$ (Lemma 5
in \cite{vdB4}).
\begin{lemma}\label{Lem1}
Let $M$ be a smooth, geodesically complete Riemannian manifold
with non-negative Ricci curvature, and let $D$ be an open subset
of $M$ with boundary $\partial D$. Then for $x\in D, t>0$
\begin{equation*}\label{e44}
    u_1(x;t)\ge 1-2^{(2+m)/2}e^{-\delta(x)^2/(8t)}.
\end{equation*}
\end{lemma}

To prove \eqref{e43} we first consider the case
\begin{equation}\label{e45}
(2+m-d)/2<\alpha<2.
\end{equation}
This set of $\alpha$'s is non-empty since $d\in [m-1,m)$. By
\eqref{e8} we have that
\begin{align}\label{e46}
\lVert\psi_{\alpha} h^{1/q}\rVert_{q/(q-1)}=c^{2/q}
\left(\int_D\delta^{(2-q\alpha)/(q-1)}\right)^{(q-1)/q}.
\end{align}
Denote the left hand side of \eqref{e42} by $\omega_D(\rho)$. Then
we can write the right hand side of \eqref{e46} as
\begin{equation}\label{e46a}
c^{2/q}\left(\int_{\R^+}
\rho^{(2-q\alpha)/(q-1)}\omega_D(d\rho)\right)^{(q-1)/q}.
\end{equation}
An integration by parts, using \eqref{e45} shows that \eqref{e46a}
is finite for
\begin{equation}\label{e47}
q<\frac{2-m+d}{\alpha-m+d}.
\end{equation}
Since $\alpha$ satisfies \eqref{e45}, we have that the right hand
side of \eqref{e47} is in $(1,2)$. We now choose $\epsilon
>0$ such that
\begin{equation}\label{e48}
\frac{2-m+d}{\alpha-m+d}(1+\epsilon)^{-1} \in (1,2),
\end{equation}
and choose $q$ equal to the left hand side of \eqref{e48}. By
Lemma \ref{Lem1} and \eqref{e42} we have that for $t\rightarrow 0$
\begin{align}\label{e49}
\lVert 1-u_1(\cdot;t)\rVert_1&=\int_D(1-u(x;t))dx\nonumber \\&\le
2^{(m+2)/2}\int_De^{-\delta^2/(8t)}\nonumber \\
&\le2^{(m+2)/2}\int_{\R^+}e^{-\rho^2/(8t)}\omega_D(d\rho)\nonumber
\\ &=2^{(m+2)/2}e^{-\rho_D^2/(8t)}|D|+2^{(m-2)/2}\kappa_Dt^{-1}\int_0^{\rho_D}\rho^{m-d+1}
e^{-\rho^2/(8t)}d\rho\nonumber \\ &=O(t^{(m-d)/2}).
\end{align}
By Theorem \ref{The1} and \eqref{e46}-\eqref{e49} we find that for
all $\alpha$ satisfying \eqref{e45} and all $\epsilon >0$
satisfying \eqref{e48}
\begin{equation}\label{e50}
Q_{\psi_{\alpha},\psi_{\alpha}}(t)=O(t^{-\epsilon(\alpha-m+d)+(m-d-2\alpha)/2}),\
t\rightarrow 0 .
\end{equation}
We conclude that \eqref{e43} holds for all $\alpha=\beta$
satisfying \eqref{e45}, and all $\epsilon>0.$

Next consider the case
\begin{equation}\label{e51}
1<\alpha<(2+m-d)/2.
\end{equation}
This set of $\alpha$'s is again non-empty since $d\in [m-1,m)$. By
\eqref{e42} we have that
\begin{equation}\label{e52}\lVert\psi^r\rVert_{q}=\left(\int_{\R+} \omega_D(d\rho) \rho^{-\alpha
rq}\right)^{1/q}< \infty \end{equation} for
\begin{equation}\label{e53}
\alpha rq<m-d,
\end{equation}
and \begin{equation}\label{e54} \lVert\psi^{2-r}h^{1/q}\rVert
_{q/(q-1)}=\left(\int_{\R^+} \omega_D(d\rho) \rho^{(2-\alpha
(2-r)q)/(q-1)}\right)^{(q-1)/q}< \infty \end{equation} for
\begin{equation}\label{e55}
\frac{\alpha q(2-r)-2}{q-1}<m-d.
\end{equation}
The optimal choice for $r$ is henceforth given by
\begin{equation}\label{e56}
r=2(\alpha q-1)\alpha^{-1}q^{-2}.
\end{equation}
By \eqref{e51} we also have that $\alpha >1$. Hence $r \in (0,2).$
The requirements under \eqref{e53} and \eqref{e55} become with
this choice of $r$ that
\begin{equation}\label{e57}
q<2(2\alpha+d-m)^{-1}.
\end{equation}
Since $\alpha$ satisfies \eqref{e51}, the right hand side of
\eqref{e57} is in $(1,2)$. We now choose $\epsilon >0$ such that
\begin{equation}\label{e58}
2((2\alpha+d-m)(1+2\epsilon))^{-1}\in (1,2),
\end{equation}
and choose $q$ equal to the left hand side of \eqref{e58}. By
Theorem \ref{The3} and \eqref{e52}-\eqref{e57} we find that for
all $\alpha$ satisfying \eqref{e51}, and all $\epsilon >0$
satisfying \eqref{e58}
\begin{equation}\label{e59}
Q_{\psi_{\alpha},\psi_{\alpha}}(t)=O(t^{-\epsilon(2\alpha-m+d)+(m-d-2\alpha)/2}),\
t\rightarrow 0.
\end{equation}
We conclude that \eqref{e43} holds for all $\alpha=\beta$
satisfying \eqref{e51}, and all $\epsilon>0$.

To prove \eqref{e43} for the limiting case
$\alpha=\beta=(2+m-d)/2:=\alpha_c$ we note that $
Q_{\psi,\phi}(t)$ is monotone on the positive cone of non-negative
and measurable $\psi$ and $\phi$. Let $\alpha=\alpha_c+\epsilon$
where $\epsilon$ is such that $\alpha\in(\alpha_c,2)$. Since
\begin{equation*}\label{e60}
\psi_{\alpha_c}\le \rho_D^{\alpha - \alpha_c}\psi_{\alpha}.
\end{equation*}
we have by \eqref{e50} that
\begin{align}\label{e61}
Q_{\psi_{\alpha_c},\psi_{\alpha_c}}(t)&\le \rho_D^{2(\alpha -
\alpha_c)}Q_{\psi_{\alpha},\psi_{\alpha}}(t)\nonumber
\\&\le\rho_D^{2(\alpha -
\alpha_c)}O(t^{-\epsilon(\alpha-m+d)+(m-d-2\alpha)/2})\nonumber
\\ &= O (t^{-\epsilon(2+\epsilon+(d-m)/2)+(m-d-2\alpha_c)/2}).
\end{align}
We conclude that \eqref{e43} holds for $\alpha=\beta=\alpha_c$,
and all $\epsilon>0$.
\end{proof}

\section{The special case calculation for a ball in $\R^3$ \label{sec4}}

In this section we show by means of an example that the upper
bound obtained in Theorem \ref{The5} is close to being sharp for
$\alpha<2, \beta<2, \alpha+\beta>3$.
\begin{theorem}\label{The6} Let $B_a=\{x\in \R^3:|x|<a\}$.
If $\alpha<2, \beta<2, \alpha+\beta>3, J\in \N$ then there exist
coefficients $b_0, b_1,\cdots$ depending on $\alpha$ and $\beta$
only such that for $t \rightarrow 0$
\begin{align}\label{e65}
Q_{\psi_{\alpha},\psi_{\beta}}&(t)=4\pi c_{\alpha,
\beta}a^2t^{(1-\alpha
-\beta)/2}-4\pi(c_{\alpha-1,\beta}+c_{\alpha,\beta-1})at^{(2-\alpha -\beta)/2} \nonumber \\
&+4\pi c_{\alpha-1,\beta-1}t^{(3-\alpha
-\beta)/2}+\sum_{j=0}^Jb_ja^{3-j-\alpha-\beta}t^{j/2}+O(t^{(J+1)/2}),
\end{align}
where
\begin{align}\label{e66}
c_{\alpha, \beta}=&2^{-\alpha -\beta}\pi^{-1/2}\Gamma((2-\alpha
-\beta)/2)\nonumber \\ & \times\int_0^1
(\rho^{-\alpha}+\rho^{-\beta})((1-\rho)^{\alpha+\beta-2}-(1+\rho)^{\alpha+\beta-2})d\rho,
\end{align}
and
\begin{align}\label{e66a}
b_0&=-8\pi((\alpha+\beta-1)(\alpha+\beta-2)(\alpha+\beta-3))^{-1},\nonumber
\\
b_1&=0, \nonumber \\
b_2&=8\pi \alpha \beta
 ((\alpha+\beta+1)(\alpha+\beta)(\alpha+\beta-1))^{-1}, \nonumber
 \\
 b_3&=0.
\end{align}
\end{theorem}
We see that the leading term in \eqref{e65} jibes with \eqref{e43}
since \eqref{e8} holds for some $c\ge 2$, and \eqref{e42} holds
with $d=m-1$.

Theorem \ref{The6} suggests that for any precompact $D$ with
smooth $\partial D$ in $M$, and for $\alpha<2, \beta<2,
\alpha+\beta>3$ and $t\rightarrow 0$
\begin{align}\label{e66b}
Q_{\psi_{\alpha},\psi_{\beta}}(t)=& c_{\alpha,\beta}\int_{\partial
D}t^{(1-\alpha-\beta)/2}-2^{-1}(c_{\alpha-1,\beta}+c_{\alpha,\beta-1})\int_{\partial
D}L_{gg}t^{(2-\alpha-\beta)/2}\nonumber \\ &+\int_{\partial
D}(c_1L_{gg}L_{hh}+c_2L_{gh}L_{gh})t^{(3-\alpha-\beta)/2}+O(1),
\end{align}
where $c_1$ and $c_2$ are constants depending on $\alpha$ and
$\beta$ only, and which satisfy
\begin{equation*}\label{e66c}
4c_1+2c_2=c_{\alpha-1,\beta-1},
\end{equation*}
and where $L_{gg}$ is the trace of the second fundamental form on
the boundary of $\partial D$ oriented by an inward unit vector
field. Since $\int_{\partial B_a}1=4\pi a^2$, $\int_{\partial B_a}
L_{gg}=8\pi a$ and $\int_{\partial
B_a}(c_1L_{gg}L_{hh}+c_2L_{gh}L_{gh})=16\pi c_1+8\pi c_2$, we see
that \eqref{e66b} holds for the ball in $\R^3$.

The proof of Theorem \ref{The6} rests on the following result
(pp.237, 367-368 in \cite{CJ}).
\begin{lemma}\label{Lem2} Let $B_a$ as in Theorem \ref{The6}, and let the initial datum be radially
symmetric i.e. $\psi_1(x)=f(r)$, where $r=|x|$. Then the solution
of \eqref{e1}, \eqref{e2}, \eqref{e3} is given by
\begin{equation*}\label{e67}
u(x;t)=(4\pi tr^2)^{-1/2}\int_0^a r'f(r')\sum_{n \in
\Z}(e^{-(2na-r+r')^2/(4t)}-e^{-(2na+r+r')^2/(4t)})dr'.
\end{equation*}
\end{lemma}

To prove Theorem \ref{The6} we have by Lemma \ref{Lem2} that
\begin{align}\label{e68}
Q_{\psi_{\alpha},\psi_{\beta}}(t)=&(4\pi/t)^{1/2}\iint_{S_a}rr'(a-r)^{-\alpha}(a-r')^{-\beta} \nonumber \\
&\times\sum_{n \in
\Z}(e^{-(2na-r+r')^2/(4t)}-e^{-(2na+r+r')^2/(4t)})drdr',
\end{align}
where $S_a=[0,a]\times[0,a]$. Substitution of $a-r=p$ and $a-r'=q$
in \eqref{e68} gives that
\begin{equation*}\label{e69}
Q_{\psi_{\alpha},\psi_{\beta}}(t)=A_0+A_1+A_2+B,
\end{equation*}
where
\begin{equation*}\label{e70}
A_0=(4\pi
/t)^{1/2}a^2\iint_{S_a}p^{-\alpha}q^{-\beta}(e^{-(p-q)^2/(4t)}-e^{-(p+q)^2/(4t)})dpdq,
\end{equation*}
\begin{equation*}\label{e71}
A_1=-(4\pi/
t)^{1/2}a\iint_{S_a}(p+q)p^{-\alpha}q^{-\beta}(e^{-(p-q)^2/(4t)}-e^{-(p+q)^2/(4t)})dpdq,
\end{equation*}
\begin{equation*}\label{e72}
A_2=(4\pi/
t)^{1/2}\iint_{S_a}p^{1-\alpha}q^{1-\beta}(e^{-(p-q)^2/(4t)}-e^{-(p+q)^2/(4t)})dpdq,
\end{equation*}
and
\begin{align}\label{e73}
B=(4\pi/ t)^{1/2}&\iint_{S_a}(a-p)(a-q)p^{-\alpha}q^{-\beta}
\sum_{n\ge 1}(e^{-(2na+p-q)^2/(4t)}\nonumber
\\ &+e^{-(2na+q-p)^2/(4t)}-e^{-(2na+q+p)^2/(4t)}-e^{-(2na-q-p)^2/(4t)})dpdq.
\end{align}
We have the following.
\begin{lemma}\label{lem3}If $1<\alpha<2, 1<\beta<2$ then for
$t\rightarrow 0$
\begin{equation}\label{e74}
B=-8\pi^{1/2}3^{-1}a^{-\alpha-\beta}t^{3/2}+O(t^2).
\end{equation}
\end{lemma}
\begin{proof}
The integrand in \eqref{e73} can be rewritten as
\begin{align}\label{e75}
&(a-p)(a-q)p^{-\alpha}q^{-\beta}\sum_{n\ge
1}e^{-(2na-p-q)^2/(4t)}\nonumber \\
&\times((e^{(p-2na)q/t}+e^{(q-2na)p/t})(1-e^{-pq/t})-(1-e^{-2pna/t})(1-e^{-2qna/t})).
\end{align}
The contribution from the terms with $n\ge2$ in \eqref{e75} is
bounded in absolute value by
\begin{equation*}\label{e76}
2a^2p^{1-\alpha}q^{1-\beta}t^{-1}\sum_{n\ge2}e^{-a^2(n-1)^2/t}(1+2n^2a^2t^{-1}).
\end{equation*}
After integrating with respect to $p$ and $q$ we see that this
term contributes at most $O(e^{-a^2/(2t)})$ to $B$. Next we will
show that the main contribution from the term with $n=1$ in
\eqref{e75} comes from a neighbourhood of the point $(p,q)=(a,a)$.
Let
\begin{equation*}\label{e77}
C_1(a)=\{(p,q)\in \R^2: a/3<p<a,a/3<q<a\},
\end{equation*}
and
\begin{equation*}\label{e78}
C_2(a)=S_a\setminus C_1(a).
\end{equation*}
On $C_2(a)$ we have that $2a-p-q\ge 2a/3$. Hence the term with
$n=1$ in \eqref{e75} is bounded on $C_2(a)$ in absolute value by
\begin{equation}\label{e79}
2(a-p)(a-q)p^{1-\alpha}q^{1-\beta}t^{-1}e^{-a^2/(9t)}(1+2a^2t^{-1}).
\end{equation}
Integrating \eqref{e79} over $C_2(a)$ gives a contribution which
is bounded by \\ $O(e^{-a^2/(18t)})$. In order to calculate the
contribution from the term with $n=1$ on $C_1(a)$ we use the
expression under \eqref{e73} instead. First we note that
$2a+p-q\ge 2a/3, 2a+q-p\ge 2a/3, 2a+p+q \ge 8a/3$. Hence the first
three terms in the summand of \eqref{e73} with $n=1$ give after
integration over $C_1(a)$ a contribution $O(e^{-a^2/(18t)})$.
Putting all this together gives that
\begin{align*}%\label{e80}
B=&-(4\pi/t)^{1/2}\iint_{C_1(a)}
(a-p)(a-q)p^{-\alpha}q^{-\beta}\nonumber
\\ & \times e^{-(2a-q-p)^2/(4t)}dpdq+O(e^{-a^2/(18t)}).
\end{align*}
Noting that
\begin{equation}\label{e81}
p^{-\alpha}q^{-\beta}=a^{-\alpha-\beta}+O(a-p)+O(a-q)
\end{equation}
uniformly in $p$ and $q$ yields after a change of variables that
\begin{align*}%\label{e82}
B=&-(4\pi/t)^{1/2}a^{-\alpha-\beta}\iint_{S_{a/3}}pqe^{-(p+q)^2/(4t)}\nonumber \\
&\times (1+O(p)+O(q))dpdq+O(e^{-a^2/(18t)}),
\end{align*}
which agrees with the right hand side of \eqref{e74}.
\end{proof}
By taking higher order terms of the form $(a-p)^{n_1}(a-q)^{n_2}$
in \eqref{e81} into account one can determine the coefficient
$t^{(j+3)/2},j=0,1,2,\cdots$ in the expansion of $B$.

To complete the proof of Theorem \ref{The6} we rewrite $A_0,A_1$
and $A_2$ respectively as follows.
\begin{align}\label{e83}
A_0=&(4\pi/t)^{1/2}a^2\left(\int_0^adp\int_0^pdq+\int_0^adq\int_0^qdp\right)\nonumber
\\ & \times
p^{-\alpha}q^{-\beta}(e^{-(p-q)^2/(4t)}-e^{-(p+q)^2/(4t)})\nonumber
\\
=&(4\pi/t)^{1/2}a^2\int_0^ap^{1-\alpha-\beta}dp\int_0^1(\rho^{-\alpha}+\rho^{-\beta})\nonumber
\\
&\times(e^{-p^2(1-\rho)^2/(4t)}-e^{-p^2(1+\rho)^2/(4t)})d\rho\nonumber
\\ =& 4\pi a^2c_{\alpha,\beta}t^{(1-\alpha-\beta)/2}\nonumber \\
&-(4\pi/t)^{1/2}a^2\int_a^{\infty}p^{1-\alpha-\beta}dp\int_0^1(\rho^{-\alpha}+\rho^{-\beta})\nonumber
\\ &\times(e^{-p^2(1-\rho)^2/(4t)}-e^{-p^2(1+\rho)^2/(4t)})d\rho,
\end{align}
\begin{align}\label{e84}
A_1=&-4\pi
a(c_{\alpha-1,\beta}+c_{\alpha,\beta-1})t^{(2-\alpha-\beta)/2}
+(4\pi/t)^{1/2}a\int_a^{\infty}p^{2-\alpha-\beta}dp\nonumber
\\ &\times \int_0^1d(\rho^{1-\alpha}+\rho^{-\alpha}+\rho^{1-\beta}+\rho^{-\beta})
(e^{-p^2(1-\rho)^2/(4t)}-e^{-p^2(1+\rho)^2/(4t)})d\rho,
\end{align}
and
\begin{align}\label{e85}
A_2&=4\pi
c_{\alpha-1,\beta-1}t^{(3-\alpha-\beta)/2}-(4\pi/t)^{1/2}\int_a^{\infty}p^{3-\alpha-\beta}dp\nonumber
\\ &\times \int_0^1d(\rho^{1-\alpha}+\rho^{1-\beta})
(e^{-p^2(1-\rho)^2/(4t)}-e^{-p^2(1+\rho)^2/(4t)})d\rho.
\end{align}
The terms to be evaluated in \eqref{e83}, \eqref{e84} and
\eqref{e85} are all of the form
\begin{equation}\label{e86}
(4\pi/t)^{1/2}a^{2-j}\int_a^{\infty}p^{1+j-\alpha-\beta}dp\int_0^1
\rho^{-
\gamma}(e^{-p^2(1-\rho)^2/(4t)}-e^{-p^2(1+\rho)^2/(4t)})d\rho,
\end{equation}
where $j=0,1,2$ respectively. Following arguments similar to the
proof of Lemma \ref{lem3} we see that the contribution of the
integral with respect to $\rho \in [0,1/2)$ in \eqref{e86} is at
most $O(e^{-a^2/(18t)})$. Furthermore
\begin{equation}\label{e86a}
(4\pi/t)^{1/2}a^{2-j}\int_a^{\infty}p^{1+j-\alpha-\beta}dp\int_{1/2}^1
\rho^{- \gamma}e^{-p^2(1+\rho)^2/(4t)}d\rho=O(e^{-a^2/(18t)}).
\end{equation}
Hence the expression under \eqref{e86} equals
\begin{equation}\label{e86b}
(4\pi/t)^{1/2}a^{2-j}\int_a^{\infty}p^{1+j-\alpha-\beta}dp\int_{1/2}^1
\rho^{- \gamma}e^{-p^2(1-\rho)^2/(4t)}d\rho+O(e^{-a^2/(18t)}).
\end{equation}
Expanding $\rho^{-\gamma}$ about $\rho=1$ we obtain that
\begin{align}\label{e87}
&|\rho^{-\gamma}-1-\gamma(1-\rho)-2^{-1}\gamma
(\gamma+1)(1-\rho)^2\nonumber \\ &-6^{-1}\gamma
(\gamma+1)(\gamma+2)(1-\rho)^3|\le C(1-\rho)^4,\ 0\le \rho \le
1/2,
\end{align}
where $C$ depends on $\gamma$ only. By \eqref{e87} and
\eqref{e86b} we obtain that \eqref{e86} is equal to
\begin{align}\label{e88}
&2\pi(\alpha+\beta-j-1)^{-1}a^{3-\alpha-\beta}+4\pi
^{1/2}\gamma(\alpha+\beta-j)^{-1}a^{2-\alpha-\beta}t^{1/2}\nonumber \\
& +2\pi
\gamma(\gamma+1)(\alpha+\beta-j+1)^{-1}a^{1-\alpha-\beta}t\nonumber
\\ &
+8\pi^{1/2}3^{-1}\gamma(\gamma+1)(\gamma+2)(\alpha+\beta-j+2)^{-1}a^{-\alpha-\beta}t^{3/2}+O(t^2).
\end{align}
It remains to compute the coefficients $b_0,b_1$ and $b_2$ in
Theorem \ref{The6}. Altogether there are eight terms which
contribute to the terms in \eqref{e88}:
$$\begin{array}{lllll}
j=0,& \gamma=\alpha,& \gamma =\beta\\
 j=1,& \gamma=\alpha-1,&
\gamma=\beta -1,& \gamma=\alpha,& \gamma=\beta \\
 j=2,&
\gamma=\alpha-1,& \gamma=\beta-1\,. \end{array}$$
 Summing these
eight terms yield the expressions for $b_0,b_1$ and $b_2$ under
\eqref{e66a}. To calculate $b_3$ we have that the above eight
$\gamma(\gamma+1)(\gamma+2)$ terms in \eqref{e88} cancel the
contribution from \eqref{e74}. This completes the proof of Theorem
\ref{The6}.

\end{document}